\documentclass[11pt]{article}
\usepackage[T1]{fontenc}
\usepackage{lmodern,amsmath,amssymb,amsthm,mathtools,booktabs,longtable,array,graphicx,verbatim,xurl,fvextra}
\usepackage[margin=1in]{geometry}
\usepackage{microtype}
\usepackage{mathrsfs,framed}
\usepackage{etoolbox}
\makeatletter
\patchcmd{\l@section}{1.5em}{2.2em}{}{
  \PackageError{manuscript}{Could not widen contents section numbers}{}
}
\renewcommand{\l@subsection}{\@dottedtocline{2}{1.5em}{3.8em}}
\makeatother
\usepackage[pdfusetitle,colorlinks=true,linkcolor=black,citecolor=blue,urlcolor=blue]{hyperref}
\usepackage{bookmark}
\apptocmd{\thebibliography}{\addcontentsline{toc}{section}{\refname}}{}{%
  \PackageError{manuscript}{Could not index the references}{}}
\newtheorem{theorem}{Theorem}[section]
\newtheorem{proposition}[theorem]{Proposition}
\newtheorem{lemma}[theorem]{Lemma}
\newtheorem{corollary}[theorem]{Corollary}
\theoremstyle{definition}
\theoremstyle{remark}
\theoremstyle{definition}
\newtheorem{mathematicalalgorithm}{Algorithm}[section]
\newcommand{\R}{\mathbb R}\newcommand{\Z}{\mathbb Z}
\newcommand{\tr}{\operatorname{tr}}\newcommand{\diag}{\operatorname{diag}}
\newcommand{\covol}{\operatorname{covol}}
\newcommand{\Sym}{\operatorname{Sym}}

\DeclarePairedDelimiter{\abs}{\lvert}{\rvert}

\theoremstyle{definition}\newtheorem{example}[theorem]{Example}
\newenvironment{planesalgorithmblock}
 {\par\addvspace{\medskipamount}\noindent\begin{minipage}{\linewidth}
  \hrule\smallskip}
 {\smallskip\hrule\end{minipage}\par\addvspace{\medskipamount}}
\newenvironment{algorithmsteps}
 {\begin{enumerate}\setlength{\itemsep}{1pt}\setlength{\parsep}{0pt}\setlength{\parskip}{0pt}}
 {\end{enumerate}}
\title{Local maxima of the rank-two {Rankin} invariant in dimension four}
\author{Avi Mehra\\University of California, Berkeley}
\date{}
\begin{document}
\maketitle

\begin{abstract}
The rank-four Euclidean lattices that locally maximize the least squared
area of a rank-two sublattice at fixed covolume are $A_4$, $A_4^*$ and $D_4$,
up to similarity. Among any nine minimum
planes of a rank-four lattice, three share a line and span the ambient space, or lie in a
three-dimensional subspace and have zero common intersection. The bound
is sharp, and the eight-plane configuration avoiding both patterns is
unique up to equivalence and occurs in $D_4$. An index bound for six
independent minimum bivectors reduces the classification to an exact
finite search. For $D_5$, ten independent minimum bivectors generate a
sublattice of the exterior square of index at most four. The selections
attaining four form eighteen orbits under lattice automorphisms, only three
of which contain six independent minimum bivectors from a copy of $D_4$. In $D_n$ the
indices of independent minimum bivectors are unbounded as $n$ grows.
\end{abstract}

\section*{Introduction}
\addcontentsline{toc}{section}{Introduction}

Which rank-four lattices locally maximize the least squared area of a
lattice plane at fixed covolume? Rankin proved that $D_4$ is the unique
global maximizer up to similarity, with normalized value $3/2$
\cite[Theorem~5]{Rankin}. Coulangeon's criterion characterizes local
maxima by $2$-perfection and $2$-eutaxy of their minimum planes and gives
$A_4$ and $D_4$ as examples \cite[Theorems~3.2.3 and~5.1.1]{Coulangeon}.
Here $2$-eutaxy, abbreviated to eutaxy, means that the orthogonal
projectors $\Pi_P$ onto the minimum planes satisfy
$\sum_P a_P\Pi_P=I$ with $a_P>0$ for every minimum plane $P$.
Watanabe's table of $2$-perfect examples includes these lattices,
$A_4^*$, and the non-eutactic pair $W_4,W_4^*$
\cite[Appendix~A6, Table~1]{Watanabe2025}.
Theorem~\ref{thm:rankfour} shows that the local maxima are exactly $A_4$, $A_4^*$ and $D_4$, up to similarity.

Our reduction uses two sharp incidence bounds. In rank four, any nine minimum planes
contain a primal or dual star, as do any eight whose bivectors span
the exterior square. Eight avoiding both stars determine the lattice
up to similarity and the selection up to automorphism
(Theorem~\ref{thm:stars}). The classification in \S\ref{sec:rankfour} also
uses the fact that six independent minimum bivectors generate a sublattice of
index at most two in the exterior square. Section~\ref{sec:extensions} shows
how this bound fails in higher rank: the maximum index is already four in
$D_5$, and the indices attained in $D_n$ are unbounded as $n$ grows.

At charge three the classical exterior-square action of
$\operatorname{SL}_4(\mathbb R)$ on $\mathbb R^{3,3}$
\cite[Appendix~A]{AchmedZadeTFolds} realizes Bavard's even split family
$V^{\mathrm{II}}_{3,3}$ as exterior squares of covolume-one rank-four
lattices \cite[\S4.9, (61)]{BavardIsodual}.
For the Narain spectral gap, Chu and Jung identify twice the gap with the
least squared plane area and hence its local maxima with the images of
the $2$-extreme lattices \cite[Corollary~9 and Theorem~11]{ChuJung}.
Theorem~\ref{thm:rankfour} exhausts these lattices. Since dual lattices
have equivalent exterior squares, the three similarity classes give
exactly two Narain local maxima up to $O(3)\times O(3)$, represented
by the normalized exterior squares of $D_4$ and $A_4$, with gaps $3/4$
and $3/(2\sqrt5)$, respectively \cite{NarainExtremaCandidate}.
\clearpage
\begingroup\small\tableofcontents\endgroup

\section{Rank-four classification}\label{sec:rankfour}

For independent
vectors $u,v$ of a Euclidean lattice $L$, the sublattice $P=\Z u+\Z v$ has
squared covolume
\[
 \covol(P)^2=\|u\wedge v\|^2
 =\|u\|^2\|v\|^2-\langle u,v\rangle^2.
\]
It is \emph{primitive} if
$P=L\cap\operatorname{span}_{\R}P$.
Passing to this full intersection can only decrease the area, so the
minimum can be taken over primitive sublattices. Set
\[
 \mu_2(L)=\min_{P\text{ primitive rank }2}\covol(P)^2,
 \]
and, in rank four, write $h_2(L)=\mu_2(L)/\covol(L)$.
This ratio is scale invariant: replacing $L$ by $tL$, for $t>0$, multiplies
both the squared plane area and the four-dimensional covolume by $t^4$.
For a Gram matrix $G$ of $L$, this denominator is $\sqrt{\det G}$.
We call the real span of a primitive sublattice attaining $\mu_2(L)$ a
\emph{minimum plane}. For a basis $u,v$ of such a sublattice, the bivectors
$\pm(u\wedge v)$ are its \emph{minimum bivectors}, one for each orientation
of the plane.

\begin{theorem}[Classification of rank-four 2-extreme lattices]\label{thm:rankfour}
On the space of rank-four Euclidean lattices up to similarity, the local
maxima of $h_2$ are precisely the classes of $A_4$, $A_4^*$ and $D_4$,
with values $3/\sqrt5$, $3/\sqrt5$ and $3/2$, respectively.
\end{theorem}

The lattices are $A_4=\{z\in\Z^5:\sum_i z_i=0\}$ and
$D_4=\{z\in\Z^4:\sum_i z_i\in2\Z\}$ with their inherited Euclidean
metrics, and $A_4^*$ is the Euclidean dual in the span of $A_4$.

At a local maximum there are at least ten minimum planes. If some eight
contain neither a primal nor a dual star, Theorem~\ref{thm:stars} identifies
the lattice as $D_4$ up to similarity. Otherwise every eight contain a star.
In this case the index bound for six independent minimum bivectors reduces
the classification to a finite list of configurations
(Lemma~\ref{prop:rankfour-configuration-census}).
The \hyperref[para:metric-outcome]{metric tests} leave one reference, a model for the pairing,
integral coordinate lattice and selected directions. It is realized by $A_4$. Metric uniqueness and integral exterior-square
equivalence then give $A_4$ or its dual
(Lemma~\ref{lem:unique-config-metric} and Proposition~\ref{prop:integral-exterior-equivalence}).

Throughout this argument, write $V=\covol(L)$, $m=\mu_2(L)$ and
$t=h_2(L)=m/V$. The exterior lattice $K=\bigwedge^2L$ carries the integral
wedge pairing $Q$, normalized by a chosen primitive volume element;
$s(x)=Q(x,x)/2$ is its quadratic value. We also write $Q$ for its matrix in a chosen
integral basis.

\subsection{The Rankin extremality criterion}
\begin{proposition}[Coulangeon's criterion for local maximality]\label{prop:rankin-direct}
Let $L$ have rank $n\ge3$, and let $\Pi_P$ denote the orthogonal projector
onto each minimum plane $P$ of $L$. Then $L$ is a local maximum of
$h_2(L)=\mu_2(L)/\covol(L)^{4/n}$ if and only if the $\Pi_P$ span
the space $\Sym_n$ of symmetric matrices and
$\sum_P a_P\Pi_P=I_n$ for some $a_P>0$.
In particular, a local maximum has at least $n(n+1)/2$ minimum planes,
hence at least ten in rank four.
\end{proposition}
\begin{proof}
The equivalence is Coulangeon's perfection--eutaxy criterion
\cite[Theorem~3.2.3]{Coulangeon} in the projector formulation
of~\cite[Section~5]{BachocCoulangeonNebe2002}. The count follows
from $\dim\Sym_n=n(n+1)/2$.
\end{proof}

\label{sec:two-perfection-comparison}
To prove that the minimum bivectors span the exterior square, we first
express Proposition~\ref{prop:rankin-direct} in bivector coordinates.

The Hodge star on bivectors is characterized by $\langle a,b\rangle=VQ(a,*b)$.
Write a bivector as $(x,y)$ in orthonormal bases of the $+1$ and $-1$
eigenspaces of $*$, respectively. Each eigenspace has dimension three.

\begin{corollary}[Coulangeon's criterion in Hodge coordinates]\label{prop:local-criterion}
A rank-four lattice is a local maximum of $h_2$ if and only if the tensors
$xy^T$ of its minimum bivectors $(x,y)$ in orthonormal Hodge coordinates
span $\R^{3\times3}$ and have a strictly positive linear combination equal
to zero.
\end{corollary}
\begin{proof}
For a trace-zero symmetric matrix $S\in\Sym_4^0$, the infinitesimal action of $e^{tS}$ on bivectors
is symmetric and anticommutes with $*$ since $e^{tS}$ preserves the wedge
pairing. Write $A(S)$ for its upper off-diagonal Hodge block; the diagonal blocks vanish.
The map $S\mapsto A(S)$ is an isomorphism: both spaces have dimension nine,
and if $A(S)=0$ then all pairwise sums of eigenvalues of $S$ vanish, so $S=0$.

For a minimum plane $W$ with bivector $(x,y)$ of squared norm $m$,
differentiating its squared norm along $e^{tS}$ in Hodge coordinates
and by the Gram determinant on $W$ gives, with the Frobenius inner
product $\langle A,B\rangle_F=\tr(A^TB)$,
\[
4\langle xy^T,A(S)\rangle_F
=2m\langle\Pi_W-I_4/2,S\rangle_F.
\]
This isomorphism identifies the two spanning conditions. It also identifies
positive relations among the tensors with positive relations among the
centered projectors. A positive relation among the latter
rescales to $\sum_W a_W\Pi_W=I_4$ and conversely. With $I_4$ in the
span, the centered projectors span $\Sym_4^0$ exactly when the
projectors span $\Sym_4$. Proposition~\ref{prop:rankin-direct} applies.
\end{proof}

\subsection{Stars among minimum planes}\label{subsec:star-incidence}
Among any nine minimum planes of a rank-four lattice, three form one of
the following two configurations (Theorem~\ref{thm:stars}).
A \emph{primal star} is a triple of planes sharing a line whose total
span is four-dimensional. A \emph{dual star} is a triple whose total span
is three-dimensional and whose common intersection is zero. A selection is star-free if it contains neither triple.
Annihilators exchange the two patterns and preserve minimality.
Indeed, extend a basis of a primitive rank-$k$ sublattice to a basis
of the rank-$n$ lattice $L$. If $G$ is its Gram matrix and $A$ its
leading $k\times k$ block, the complementary dual basis generates
the annihilator with Gram determinant $\det A/\det G$. Taking minima
on both sides gives
\begin{equation}\label{eq:planes-primitive-duality}
 \mu_{n-k}(L^*)=\mu_k(L)/V^2,\qquad
 h_{n-k}(L^*)=h_k(L),\quad h_k(L)=\mu_k(L)/V^{2k/n},
\end{equation}
where $\mu_k$ is the least squared
covolume of a primitive rank-$k$ sublattice. Duality is a smooth
involution on shapes, so it preserves strict local maximality as well.
In projective three-space, a \emph{star} consists of all lines through a
point, and a \emph{dual star} consists of all lines in a plane. Our primal
stars correspond to three noncoplanar lines in a star, and our dual stars
to three nonconcurrent lines in a dual star. Both conditions exclude
triples from a pencil, the lines in one plane through one point.\par

\begin{figure}[!h]
\centering
\includegraphics{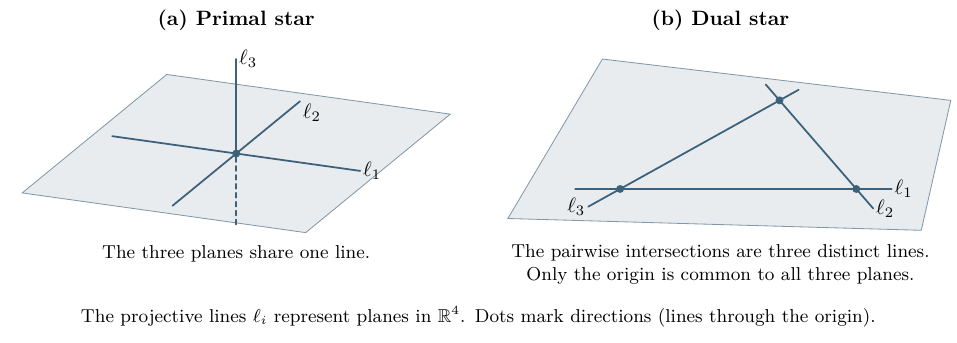}
\caption{Projectively, the three lines on the left meet at a point and
are not coplanar; those on the right lie in a plane and are not concurrent.
In the underlying vector space the three planes on the right intersect
only at the origin.}
\label{fig:minimum-planes-stars}
\end{figure}

\begin{lemma}\label{lem:independent-incident-star}
Let $P_1,P_2,P_3$ be distinct two-dimensional subspaces of a
four-dimensional real vector space $V$. For each $i$, choose a nonzero
bivector $q_i\in\bigwedge^2P_i$. The planes form a primal or dual star if and only if the $q_i$ are
linearly independent and $q_i\wedge q_j=0$ for all $i\ne j$.
\end{lemma}
\begin{proof}
The wedge products vanish exactly when the planes meet pairwise in lines.
If they share a line $\R u$, the map $v\mapsto u\wedge v$ identifies
$V/\R u$ with $u\wedge V$. Thus the $q_i$ are independent exactly when
the three planes span $V$, the primal star condition.
Otherwise $P_3\cap P_1$ and $P_3\cap P_2$ are distinct lines spanning
$P_3$, so all three planes lie in $W=P_1+P_2$. A volume form on $W$
identifies the $q_i$ with covectors having kernels $P_i$. Their common
kernel is zero, so the covectors are independent and the planes form
a dual star.
\end{proof}

On the Klein quadric, the three points $[q_i]$ span an $\alpha$-plane
in the primal case and a $\beta$-plane in the dual case
\cite[\S3.1]{RudnevSelig}.

Four planes form a \emph{transverse quartet} if each pair intersects
only at the origin. Equivalently, all six wedge pairings between their
representative bivectors are nonzero.

\begin{samepage}
\begin{theorem}[Sharp star thresholds for minimum planes in rank four]\label{thm:stars}\label{thm:fresh-eight}
Every star-free selection of eight minimum planes in a rank-four lattice
is, up to an integral change of basis and uniform scale, the configuration
in the first two rows of the $D_4$ grid in Lemma~\ref{lem:d4-selected-families}. Its primitive bivectors
span a five-dimensional subspace of the exterior square.
Consequently any nine minimum planes contain a primal or dual star, as do
any eight whose bivectors span the full exterior square. Both thresholds
are sharp for selected families.
\end{theorem}
\end{samepage}
\begin{proof}
Suppose eight minimum planes are star-free. If they contain a transverse
quartet, its wedge pairings force Rankin equality by
\eqref{eq:eight-trace} and the quartet calculation below it. Otherwise
Lemma~\ref{lem:eight-finite} gives either four planes forcing the same
equality or one of the two pairings excluded by Lemma~\ref{lem:eight-exclusions}.
Thus the lattice is similar to $D_4$. Lemma~\ref{lem:d4-selected-families}
carries the selected eight to the first two rows of its grid and gives
exterior rank five. In particular any eight spanning the exterior square
contain a star.

A star-free nine would therefore lie in $D_4$, where the grid admits
no such selection. The eight planes in the first two rows contain no star, so the bound
of nine is sharp. Under the spanning hypothesis, the bound of eight is also
sharp: the seven planes in the union of the first row and the second column
contain no star, although their bivectors span the exterior square.
\end{proof}

\subsection{Enumerating minimum-plane configurations and testing their metrics}
\label{sec:rankfour-classification}

For every rank-four lattice $L$, the shortest nonzero vectors of
$K=\bigwedge^2L$ are its minimum bivectors
\cite[Corollary~2.3]{Coulangeon1997}. Indeed, in an oriented lattice
basis $e_1,\ldots,e_4$, write $e_{ij}=e_i\wedge e_j$. The only nonzero pairings of distinct basis bivectors, up to symmetry,
are $Q(e_{12},e_{34})=1$, $Q(e_{13},e_{24})=-1$ and $Q(e_{14},e_{23})=1$.
Their self-pairings vanish, so $Q$ is even, integral and unimodular of
signature $(3,3)$. A bivector $z$ is decomposable exactly when
$z\wedge z=0$, or equivalently $s(z)=0$. The Hodge star gives
\begin{equation}\label{eq:eight-hodge}
 \|z\|^2\ge V\abs{Q(z,z)}=2V\abs{s(z)}.
\end{equation}
Thus every nonsimple integral bivector has squared norm at least $2V$,
whereas Rankin's bound gives
\begin{equation}\label{eq:eight-rankin}
 m\le 3V/2 \qquad\text{\cite[Theorem~5]{Rankin}}.
\end{equation}
Every nonzero simple integral bivector is an integer multiple of the
primitive bivector of its plane, as one sees by extending a basis of the
saturated plane lattice to a basis of $L$. This proves the claim.
For minimum bivectors $x,y$, the Hodge star likewise gives
$\abs{Q(x,y)}\le\|x\|\|y\|/V=m/V\le3/2$. Integrality leaves only
$Q(x,y)\in\{0,\pm1\}$.

\subsubsection{Local maximality and six independent minima}
\label{subsec:local-six}
\begin{lemma}[Full exterior span]\label{lem:full-exterior-span}
At a local maximum the minimum bivectors span $\bigwedge^2\mathbb R^4$.
\end{lemma}
\begin{proof}
Write the minimum bivectors in orthonormal Hodge coordinates as
$v_i=(x_i,y_i)$. By Corollary~\ref{prop:local-criterion}, the tensors
$x_iy_i^T$ span the space of $3$ by $3$ matrices and there are positive
weights $\lambda_i$ with $\sum_i\lambda_i x_iy_i^T=0$.
The tensor span forces the $x_i$ and the $y_i$ to span their respective
Hodge summands. The balance equation makes the covariance block diagonal:
\[
 \sum_i\lambda_i v_iv_i^T
 =\begin{pmatrix}\sum_i\lambda_i x_ix_i^T&0\\
                  0&\sum_i\lambda_i y_iy_i^T\end{pmatrix}.
\]
Both diagonal blocks are positive definite, so the $v_i$ span the full
exterior space.
\end{proof}

Lemma~\ref{lem:full-exterior-span} is also the exterior-square case of
Berg\'e and Martinet's spanning theorem for extreme orthogonal isodual
lattices~\cite[Theorem~7.6]{BergeMartinet1995}.

In ranks $n\ge5$, a family of rational planes can satisfy the same
projector conditions while its bivectors span only a hyperplane
(Example~\ref{thm:extension-abstract-obstruction}).

\begin{lemma}[Index bound for six independent minimum bivectors]\label{lem:minbasis-index}
In any rank-four lattice, six independent minimum bivectors generate a sublattice of $K$ of
index at most two. If this index is two, the nonzero coset is represented
by half the sum of five or six of these bivectors, after changing their signs.
\end{lemma}
\begin{proof}
Scale $L$ to covolume $V=1$, so $K=\bigwedge^2L$ also has covolume one.
Let $G$ be the Euclidean Gram matrix of the six minimum bivectors $x_i$,
and let $d=[K:\sum_i\mathbb Z x_i]$. By Hadamard's inequality and
Rankin's bound~\eqref{eq:eight-rankin},
\[
 d=\sqrt{\det G}\le m^3\le(3/2)^3<4.
\]

Suppose $d=3$. Reducing coefficients modulo three and choosing signs
for the $x_i$, represent a nonzero coset by
$v=\frac13\sum_{i\in I}x_i$. Each $v-x_i$ for $i\in I$ lies in that
coset, so Watson's identity~\cite{Watson1971} gives
\[
 0\le\sum_{i\in I}(\|v-x_i\|^2-m)=(\abs{I}-6)\|v\|^2.
\]
Thus $\abs{I}=6$ and every $v-x_i$ is shortest, hence decomposable.
Expanding $Q(v-x_i,v-x_i)=0$ gives $Q(v,x_i)=Q(v,v)/2$ for every $i$.
Since the $x_i$ span, $Q$ is nondegenerate and $v\ne0$, this value is
nonzero. Evenness of $Q$ and~\eqref{eq:eight-hodge} now give
$\|v\|^2\ge\abs{Q(v,v)}\ge2$.

The equalities $\|v-x_i\|^2=\|x_i\|^2$ give
$\langle x_i,v\rangle=\|v\|^2/2$. Since $\sum_i x_i=3v$, the Gram
matrix with the chosen signs has constant row sum
$\lambda=3\|v\|^2/2\ge3$. Its trace is $6m\le9$, so its five other
eigenvalues have sum at most $9-\lambda$. By AM--GM,
\[
 9=d^2=\det G\le\lambda\left(\frac{9-\lambda}{5}\right)^5
 \le3(6/5)^5<9,
\]
where $\lambda(9-\lambda)^5$ decreases for $3\le\lambda<9$.

If $d=2$, the nonzero coset has a representative
$\frac12\sum_{i\in I}x_i$. Its support gives a zero sum in $K/2K$.
By Lemma~\ref{lem:binary-lifting} the reductions are distinct and
nonzero, and a sum of three or four would give a signed relation among
the independent $x_i$. Thus $\abs{I}$ is five or six.
\end{proof}

\subsubsection{A finite set of candidate bivectors}
\label{subsec:finite-basis}
Suppose $L$ is a local maximum and every eight minimum planes contain a star.
Lemma~\ref{lem:full-exterior-span} supplies six independent minimum bivectors.
Among all such choices take $x_1,\ldots,x_6$ whose integer span has largest index $d$ in $K$.
In an integral basis of $K$, the columns $x_i$ have absolute determinant $d$.
For a minimum bivector $x=\sum_j c_jx_j$, multilinearity gives absolute
determinant $d\abs{c_i}$ after replacing column $i$ by $x$. When $c_i\ne0$,
this is the new index, so maximality gives $\abs{c_i}\le1$.

The six pairings $b_i=Q(x_i,x)\in\{0,\pm1\}$ determine $x$ by nondegeneracy.
By Lemma~\ref{lem:minbasis-index}, $d\in\{1,2\}$. Put $Q_{ij}=Q(x_i,x_j)$,
so $Qc=b$. Decomposability gives zero diagonal. The wedge pairing on $K$
is unimodular of signature $(3,3)$,
so passage to the six columns multiplies its determinant $-1$ by $d^2$.
Thus
\[
 Q_{ii}=0,\quad Q_{ij}\in\{0,\pm1\},\quad \det Q=-d^2.
\]
The two minimum bivectors of each minimum plane give one class in the finite candidate set
\begin{equation}\label{eq:candidate-cube}
 \bigl\{c\in\R^6 : {\textstyle\sum_i} c_i x_i\in K,\ Qc\in\{0,\pm1\}^6\setminus\{0\},\
 c^TQc=0,\ \abs{c_i}\leq1\ \text{for all }i\bigr\}/\{\pm1\}.
\end{equation}
Algorithm~\ref{alg:complete-sets} uses \eqref{eq:candidate-cube} to cover
the possible complete minimum sets.

\subsubsection{Bounds for the exterior metric}
\label{subsec:metric-certificates}
The following bounds are necessary for a candidate set $S$ to be the
complete set of minimum bivectors, up to sign, at a local maximum.

Let $H$ be the exterior metric divided by $V$. For every rank-four lattice, $Q^{-1}H$ represents the Hodge star.
In an $H$-orthonormal Hodge basis, $H=I_6$ and $Q=\diag(I_3,-I_3)$,
which gives
\begin{equation}\label{eq:metric-compatibility}
 H>0,\quad HQ^{-1}H=Q,\quad H\pm Q\succeq0,\quad
 \tr(Q^{-1}H)=0.
\end{equation}
Now suppose $L$ is a local maximum and $S$ contains one minimum
bivector from each minimum plane. Normalizing the positive weights in
Corollary~\ref{prop:local-criterion} to sum to one gives
\begin{equation}\label{eq:covariance-balance}
 E=\sum_{x\in S}\lambda_xxx^T>0,\quad
 \lambda_x>0,\quad\sum_x\lambda_x=1,\quad
 HEH=QEQ,\quad x^THx=t.
\end{equation}
Indeed, write $x=(x_+,x_-)$ in the Hodge basis. The off-diagonal block of
$E$ is $\sum_x\lambda_x x_+x_-^T=0$ by Corollary~\ref{prop:local-criterion}.
Thus $E=QEQ$ in this basis, giving $HEH=QEQ$ in the marked basis.
The positive weights and Lemma~\ref{lem:full-exterior-span} give $E>0$.

The following trace inequality applies also to singular matrices. For symmetric $A$,
$\abs{A}$ denotes the positive semidefinite square root of $A^2$.
\begin{lemma}[Trace lower bound]\label{lem:trace-reference}
Let $A,Z$ be real symmetric matrices with $Z\pm A\succeq0$. Then
$\tr Z\ge\tr\abs{A}$, with equality if and only if $Z=\abs{A}$.
\end{lemma}
\begin{proof}
The bound is the symmetric case of the nuclear-norm semidefinite
formula~\cite[Eq.~(2.6)]{RechtFazelParrilo}, with both diagonal blocks equal to $Z$.
In an orthonormal eigenbasis of $A$, each $Z_{ii}\ge\abs{A_{ii}}$.
At equality every diagonal bound is sharp, so for each $i$ either
$Z-A$ or $Z+A$ has zero $i$th diagonal entry. Positivity then forces
its $i$th row to vanish. Hence $Z=\abs{A}$, which attains the bound.
\end{proof}

\begin{lemma}[Metric uniqueness for a spanning eutactic minimum set]\label{lem:unique-config-metric}
For an actual eutactic full minimum set $S$ spanning
$\mathbb R^6$, its metric $H$
is the unique minimizer of
\[
       \min_{M\pm Q\succeq0}\ \max_{x\in S}x^TMx.
\]
The minimum is $t$. Thus any feasible $U$ supplies the upper bound
$t\le\max_{x\in S}x^TUx$. If an independently checked lower bound equals
this upper bound, then necessarily $H=U$.
\end{lemma}
\begin{proof}
For $E$ in \eqref{eq:covariance-balance}, put $R=E^{1/2}$ and $A=RQR$.
If $M\pm Q\succeq0$, then $RMR\pm A\succeq0$.
Since $HEH=QEQ$, the positive semidefinite matrix $RHR$ squares to
$A^2$, so $RHR=\abs{A}$. Lemma~\ref{lem:trace-reference} now gives
\[
 \max_{x\in S}x^TMx\ge\sum_{x\in S}\lambda_x x^TMx
 =\tr(RMR)\ge\tr(RHR)=t.
\]
The metric $H$ attains $t$. Equality forces $RMR=RHR$, and $R$ is
invertible since $E>0$, so $M=H$.
\end{proof}

\subsection{Proof of the rank-four classification}
\begin{lemma}[Minimum planes in even integral lattices]\label{lem:even-determinant-three}
In a positive definite even integral lattice, every rank-two sublattice
has squared covolume at least three. If two vectors of squared norm two have
inner product $\pm1$, the minimum planes are exactly the planes
spanned by such pairs, and each pair generates its primitive plane lattice.
\end{lemma}
\begin{proof}
Chu and Jung give the determinant bound and equality case in rank
four~\cite[Lemma~10]{ChuJung}. Their binary argument applies in every rank.
A positive binary Gram determinant $ac-b^2$ is congruent to $0$ or $3$
modulo $4$, so is at least $3$. At determinant $3$, a Gauss-reduced basis
satisfies $a\le c$ and $\abs{b}\le a/2$, hence $3\ge3a^2/4$.
Since $a$ is positive and even, $a=2$; then $2c-b^2=3$ gives $c=2$ and
$\abs{b}=1$. Conversely such a pair attains determinant $3$. A proper
primitive saturation of index $d\ge2$ would have determinant $3/d^2<3$,
contradicting the same lower bound.
\end{proof}

\begin{proof}[Proof of Theorem~\ref{thm:rankfour}]
Rankin's theorem~\cite[Theorem~5]{Rankin} gives $D_4$ as the unique
global maximum up to similarity, with $h_2(D_4)=3/2$.
The roots of $A_4$ are $e_i-e_j$. Two have inner product $\pm1$
precisely when their supports meet in one coordinate, so
Lemma~\ref{lem:even-determinant-three} gives exactly the ten coordinate-triple
$A_2$ planes, each of squared area three. Since $\covol(A_4)=\sqrt5$,
$h_2(A_4)=3/\sqrt5$. Coulangeon's theorem gives strict local maximality
\cite[Theorem~5.1.1]{Coulangeon}.
Primitive duality \eqref{eq:planes-primitive-duality} gives the same
conclusion for $A_4^*$. The two lattices have twenty and ten shortest
vectors, respectively: the roots $e_i-e_j$ in $A_4$ and
$\pm(e_i-\tfrac15\mathbf1)$ in $A_4^*$.

Now let $L$ be any local maximum. If some eight minimum planes are
star-free, Theorem~\ref{thm:stars} gives $D_4$ up to similarity.
Otherwise Lemma~\ref{prop:rankfour-configuration-census} carries the
complete minimum set to one of the enumerated references by an integral
map preserving the pairing up to sign. The
\hyperref[para:metric-outcome]{metric tests} leave just one reference.
Since $h_2(A_4)=3/\sqrt5<3/2$, Theorem~\ref{thm:stars} places $A_4$ in
this second case as well. The two transports identify the integral coordinate lattices and the
selected directions, preserving the pairing up to sign. In these coordinates,
Lemma~\ref{lem:unique-config-metric} forces the normalized exterior
metrics of $L$ and $A_4$ to agree. Proposition~\ref{prop:integral-exterior-equivalence}
therefore gives $L$ up to similarity as $A_4$ or $A_4^*$.
\end{proof}

\section{Minimum bivectors in higher rank}\label{sec:extensions}

Let $L$ have rank $n$ and put $N=\binom n2$. We study the index in
$\bigwedge^2 L$ of the sublattice generated by $N$ independent minimum
bivectors of $L$. Changing the sign of any chosen bivector preserves this
sublattice.

Lemma~\ref{lem:minbasis-index} bounds the index of six independent minimum
bivectors by two in rank four. No bound independent of rank holds even at
strict local maxima of $h_2$: these indices can be arbitrarily large in
$D_n$ as $n$ grows. Local maximality of $D_n$ for $n\ge4$ follows from
\cite[Theorem~5.1.1]{Coulangeon}. We prove unboundedness by
extending selections one coordinate at a time. Already in $D_5$, the possible
indices are $1,2,3,4$, and index three cannot arise by extending a selection
from $D_4$. Example~\ref{thm:extension-abstract-obstruction} shows that
the projector conditions alone no longer force the bivectors to span.

\subsection{Coordinate extensions}\label{subsec:extension-index}

Put $V=\covol(L)$. Each singular value
of a basis matrix is a factor in exactly $n-1$ pairwise products, so
$\covol(\bigwedge^2 L)=V^{n-1}$.

Write $D_n=\{z\in\mathbb Z^n:\sum_i z_i\text{ is even}\}$ with its
usual Euclidean metric.

\begin{lemma}[Shortest exterior vectors of $D_n$]
\label{lem:extension-Dn-shell}
For $n\geq3$, the shortest nonzero vectors of $\bigwedge^2D_n$ have squared
norm three. Writing $e_{ij}=e_i\wedge e_j$, they are exactly
\[
 \pm e_{ij}\pm e_{ik}\pm e_{jk},\qquad i<j<k.
\]
The signs are independent, giving eight vectors per triple and
$8\binom n3$ in all. We call them \emph{signed coordinate triangles}
because their supports are the three edges of a triangle.
Each is the bivector of a primitive minimum plane.
\end{lemma}
\begin{proof}
For $u,v\in D_n$, the $i$th row sum of the skew matrix of $u\wedge v$
is $u_i\sum_jv_j-v_i\sum_ju_j$, which is even. This holds for every
integral sum of such wedges. Thus the odd coefficients of a bivector
in $\bigwedge^2D_n$ form the edges of a graph with even degree at every
vertex. A nonzero bivector of squared norm at most three has at most
three nonzero coefficients, each $\pm1$. Its graph must therefore be a
triangle, so its squared norm is three. This bound is attained by
\begin{equation}\label{eq:extension-Dn-wedge}
 (e_i-e_j)\wedge(e_j-e_k)=e_{ij}-e_{ik}+e_{jk}.
\end{equation}
Changing the sign of coordinate $i$, $j$ or $k$ preserves $D_n$ and
reverses the two coefficients involving that coordinate. Together with
the displayed wedge this gives four sign patterns. Reversing the wedge
factors gives their four negatives. Each wedge is shortest and hence
primitive, so its factors form a basis of the lattice in that minimum plane.
Decomposability here also follows from Coulangeon's criterion for integral
lattices \cite[Proposition~2.4]{Coulangeon1997}.
\end{proof}

A coordinate extension from $D_m$ to $D_{m+1}$ retains
$\binom m2$ independent minimum bivectors by adjoining a zero coordinate
and adds $m$ signed coordinate triangles through the new coordinate $m+1$.

\begin{theorem}[Index range for coordinate extensions]
\label{thm:extension-Dn-index}
Let $n\ge r\ge3$. Start with $\binom r2$ independent minimum bivectors
of $D_r$ of index $I$. Successive coordinate extensions to $D_n$ that preserve
independence attain exactly the indices
\[
 I2^j,\qquad j\in\mathbb Z,\quad 0\le j\le
 \lfloor(n-3)^2/4\rfloor-\lfloor(r-3)^2/4\rfloor.
\]
\end{theorem}

\begin{proof}

In one extension, order the old rows and columns first. The row matrix
is $\left(\begin{smallmatrix}A&0\\ *&B\end{smallmatrix}\right)$, where
$A$ is the old matrix and $B$ records the $m$ new triangles in the
coordinates $e_{1,m+1},\ldots,e_{m,m+1}$. The absolute determinant
is multiplied by $\abs{\det B}$ while the exterior covolume doubles from
$2^{m-1}$ to $2^m$. Thus the index multiplier is $\abs{\det B}/2$ whenever
$B$ is nonsingular.

Each row of $B$ is $e_a-\sigma e_b$ up to sign. Regard it as an edge
$ab$ of sign $\sigma$ on the $m$ vertices, allowing parallel edges.
By the signed-incidence determinant rule
\cite[Theorem~5.1(g), Lemma~8A.2]{Zaslavsky1982}, $B$ is nonsingular
exactly when each component has a single cycle and the product of the
signs on that cycle is negative. Each such component contributes a
factor of two to $\abs{\det B}$, so with $c$ components the indices satisfy
\begin{equation}\label{eq:extension-index-step}
 I_{m+1}=2^{c-1}I_m.
\end{equation}

Each component needs at least two vertices, and every
$1\le c\le\lfloor m/2\rfloor$ occurs: take $c$ disjoint pairs, join
each pair by two edges of opposite signs, and attach the remaining
vertices as leaves. Each edge $ab$ with $a<b$ lifts to the minimum
bivector $e_{ab}+e_{a,m+1}-\sigma e_{b,m+1}$ by
Lemma~\ref{lem:extension-Dn-shell}, independently of the old rows.
Thus the exponent increases by any integer from $0$ to
$\lfloor m/2\rfloor-1$ at each step. These intervals add without gaps,
and their upper endpoints telescope by
\[
 \lfloor m/2\rfloor-1
 =\lfloor(m-2)^2/4\rfloor-\lfloor(m-3)^2/4\rfloor.\qedhere
\]

\end{proof}

In $D_4$ the six triangles $\pm e_{ij}+e_{i4}+e_{j4}$,
$1\le i<j\le3$, have absolute determinant $16$ by subtracting the two
rows of each pair, hence index $16/2^3=2$. The upper endpoint of
Theorem~\ref{thm:extension-Dn-index} therefore gives index
$2^{1+\lfloor(n-3)^2/4\rfloor}$ in $D_n$.

\subsection{Possible indices and maximizing selections in \texorpdfstring{$D_5$}{D5}}
We now allow any ten independent triangles.

\begin{theorem}[Possible indices in $D_5$]
\label{thm:extension-D5-maximum}
The indices in $\bigwedge^2 D_5$ of sublattices generated by ten
independent minimum bivectors are exactly $1,2,3,4$.
\end{theorem}
\begin{proof}
By Lemma~\ref{lem:extension-Dn-shell} there are forty minimum bivectors
up to sign: ten triangle supports, each with four sign choices. If ten
of them form the rows of $M$ in standard wedge coordinates, their index
is $\abs{\det M}/\covol(\bigwedge^2D_5)=\abs{\det M}/16$. Signed coordinate
permutations act transitively on the forty directions and preserve the
index, so it suffices to enumerate ten-subsets containing one fixed row.

If $k$ chosen rows have Gram determinant $p$, any completion satisfies
$\abs{\det M}^2\le p\,3^{10-k}$: each new row contributes its squared
distance from the preceding span, at most three. Once absolute
determinant $d$ has been found, the next improvement is at least $d+16$.
The search therefore skips a branch if $p\,3^{10-k}<(d+16)^2$.
The largest absolute determinant found is $64$, so the index is at most four.
The search exhibits determinants $16,32,64$, realizing indices one, two and four.

For index three, take the all-positive triangle on each of the ten supports.
Index its row by the complementary pair. The resulting symmetric matrix $B$
has $B_{ef}=1$ exactly when the pairs $e,f$ are disjoint. Each pair has
three disjoint partners. Two distinct pairs have one common partner if
they meet and none if they are disjoint. Thus $B\mathbf1=3\mathbf1$ and
$B^2+B=2I+\mathbf1\mathbf1^T$. On $\mathbf1^\perp$ the eigenvalues belong to
$\{1,-2\}$. Since $\tr B=0$, their multiplicities are five and four.
Hence $\abs{\det B}=3\cdot2^4=48$, giving index three.
By Lemma~\ref{lem:minbasis-index} and~\eqref{eq:extension-index-step},
an extension from $D_4$ has index a power of two, so it cannot give index three.
\end{proof}

Selections are unordered, and changing the sign
of a bivector does not change the selection; equivalence means an
automorphism of the Euclidean lattice $D_5$.

\begin{theorem}[Automorphism orbits of index-four selections in $D_5$]\label{thm:extension-D5-equality}
The maximizing selections form eighteen automorphism orbits. Exactly
three consist of selections containing six independent bivectors in the
exterior square of an isometric copy of $D_4$. In each such selection
the copy is unique, and the six bivectors generate a sublattice of
index two in its exterior square. All such sets of six are isometric
as unoriented directions.
\end{theorem}
\begin{proof}
Since $D_5^*=\Z^5\cup(\tfrac12\mathbf1+\Z^5)$ and the second coset
has squared minimum $5/4$, the shortest dual vectors are $\pm e_i$.
Every automorphism of $D_5$ therefore permutes the five coordinate axes,
and every signed coordinate permutation preserves $D_5$. To count labeled maximizers, enumerate all ten-subsets as in
Theorem~\ref{thm:extension-D5-maximum}, without fixing a row and retaining
every branch that can
reach absolute determinant $64$. The $21\,312$ maximizers form eighteen
orbits under these permutations.

Every isometric copy of $D_4$ contains four orthogonal roots, the images
of $e_1\pm e_2$ and $e_3\pm e_4$. Orthogonal roots of $D_5$ have equal
or disjoint two-coordinate supports. At most two disjoint pairs fit in
five coordinates, and each supports at most two orthogonal roots. Hence
the copy spans a coordinate hyperplane and equals its intersection with
$D_5$, since both lattices have covolume two.

For each orbit representative, count the selected triangles avoiding
each coordinate. Six such triangles span the exterior square of the
corresponding coordinate $D_4$ over $\mathbb R$, since the selection is
independent. Exactly three representatives have a count of six, each for
a unique coordinate. For each such copy of $D_4$, the matrix with the six
selected bivectors as rows in its wedge coordinates has absolute
determinant $16$, giving index $16/2^3=2$.
The three sets of six unoriented bivectors lie in one orbit under
signed coordinate permutations.
\end{proof}

For $n\ge6$, Theorem~\ref{thm:extension-Dn-index} gives a lower bound
for the maximum index in $D_n$ over all selections of $\binom n2$
independent minimum bivectors.

\subsection{The projector conditions do not force exterior span}
Lemma~\ref{lem:full-exterior-span} deduces exterior span from projector
perfection and eutaxy in rank four.

\begin{example}\label{thm:extension-abstract-obstruction}
For every $n\ge5$ there is a finite family $\mathcal F$ of rational
two-planes in $\R^n$ whose orthogonal projectors span $\Sym_n$ and satisfy
$\sum_{P\in\mathcal F}\lambda_P\Pi_P=I_n$ with every $\lambda_P>0$.
Their unit bivectors span a hyperplane in $\bigwedge^2\R^n$.
No Euclidean lattice has every plane in $\mathcal F$ as a minimum plane.
\end{example}
\begin{proof}
Start with all coordinate planes except $12$ (write $ij$ for $\langle e_i,e_j\rangle$). Their bivectors already span $e_{12}^\perp$. Write $m=n-2$. Give edges joining $\{1,2\}$ to the remaining vertices weight $1/m$, and edges within the remaining set weight $(1-2/m)/(m-1)$. Since $m\ge3$ both weights are positive. Every vertex has weighted degree one, so the weighted projectors sum to $I$. They span all diagonals: a diagonal orthogonal to them has entries satisfying $d_i+d_j=0$ on a connected graph containing a triangle, forcing every entry to vanish.
 
For every pair $i<j$ choose $k\notin\{1,2,i,j\}$. Add the two planes spanned by $e_k$ and $e_i\pm e_j$. Their unit bivectors still have zero $12$ coordinate. Their projector difference is $e_i e_j^t+e_j e_i^t$ and their projector sum equals $\Pi_{ik}+\Pi_{jk}$. Give both new planes weight $\delta$ and subtract $\delta$ from the two indicated coordinate weights for each pair. A common sufficiently small $\delta>0$ keeps every weight positive after the finitely many subtractions. Thus the sum remains $I$, the projectors span all symmetric matrices, and the bivector span has exactly codimension one.

Consider the two planes spanned by $e_k$ and, respectively, $e_1$ and
$e_1+e_2$, with $k$ chosen for $(i,j)=(1,2)$. If both were minimum
planes of $L$, their intersection and sum would be rational for $L$.
Let $u$ generate $L\cap\R e_k$ and project the lattice in their sum
orthogonally along $u$. A primitive vector $w$ of the projected lattice
lifts to a plane lattice of area $\|u\|\,\|w\|$: a lift of $w$ and $u$
form its primitive basis. Thus the two indicated directions would have
shortest representatives $v,w$ of the same length making angle $\pi/4$.
But $\|v-w\|^2=(2-\sqrt2)\|v\|^2<\|v\|^2$.
\end{proof}

\appendix
\section{Finite configurations and exterior-square equivalence}\label{app:rankfour-proofs}
\subsection{Classification of star-free sets of eight minimum planes}
\label{sec:fresh-eight}

\subsubsection{From minimum planes to binary points}
\label{subsec:binary-lifting}
\begin{lemma}[Reduction modulo two and signed relations]\label{lem:binary-lifting}
Choose a primitive bivector for each minimum plane. Their reductions are
distinct nonzero isotropic points of $K/2K$. If three or four of them sum
to zero, their representatives satisfy a relation with coefficients $\pm1$.
For three planes, real dependence of the bivectors is equivalent to
dependence of their reductions.
\end{lemma}
\begin{proof}
Each minimum bivector $x_i$ is primitive and satisfies $s(x_i)=0$, so its
reduction is nonzero and isotropic. By Section~\ref{sec:rankfour-classification},
the shortest nonzero vectors of $K$ have squared norm $m$ and are simple.
If $k$ reductions sum to zero, all half-sums
$z=\tfrac12\sum_{i=1}^k\varepsilon_i x_i$ with $\varepsilon_i=\pm1$
belong to $K$. Averaging over the signs cancels the cross terms, giving
mean squared norm $km/4$. For $k=2,3$ this is less than $m$, so some
half-sum vanishes. This excludes equal residues for distinct directions
and gives the asserted relation for triples.

For $k=4$ the mean is $m$. If no half-sum vanishes, every $z$ is shortest
and hence simple. Averaging $\|z\|^2=m$ and $Q(z,z)=0$ against
$\varepsilon_i\varepsilon_j$ isolates $\langle x_i,x_j\rangle=Q(x_i,x_j)=0$
for each $i\ne j$. Four nonzero orthogonal vectors would then span a totally
isotropic four-space for a form of signature $(3,3)$, a contradiction.

Conversely, three dependent bivectors have integral coordinates, so they
satisfy a nonzero rational relation. Clear denominators and divide by
the gcd of the coefficients. At least one coefficient is then odd, so
reduction modulo two gives a nonzero relation.
\end{proof}

Three minimum planes form a star exactly when their reductions are
independent and pairwise orthogonal. Indeed, by Lemma~\ref{lem:binary-lifting}
independence lifts in both directions for triples. Since
$Q(x_i,x_j)\in\{0,\pm1\}$, a pairing vanishes modulo two exactly when
it vanishes over $\Z$. Lemma~\ref{lem:independent-incident-star} then
gives a star. Conversely, a star has three independent bivectors with
zero mutual pairings.

In exterior coordinates $(12,13,14,23,24,34)$ the reductions lie on
\[
 \mathcal Q=\{u\in\mathbb F_2^6\setminus\{0\}:
              u_1u_6+u_2u_5+u_3u_4=0\},\qquad \abs{\mathcal Q}=35.
\]
Fix the first member of each coordinate pair. The zero triple admits
eight partners and each nonzero triple admits four, giving
$8+7\cdot4-1=35$ nonzero points.

For a real symmetric matrix $A$, write $\abs{A}=(A^2)^{1/2}$;
its trace is the sum of the absolute values of the eigenvalues.

\subsubsection{Enumerating the wedge pairings of eight star-free planes}

\paragraph{Transverse quartets.}
Let $C=X^TQX$ be the pairing matrix of four minimum bivectors, taken as
columns of $X$. Since $X^THX\pm C\succeq0$, Lemma~\ref{lem:trace-reference}
and Rankin's bound give
\begin{equation}\label{eq:eight-trace}
 \tr\abs{C}\le\tr(X^THX)=4t\le6.
\end{equation}
Thus $\tr\abs{C}=6$ forces $t=3/2$ and identifies $D_4$ up to similarity
\cite[Theorem~5]{Rankin}. In particular, this happens whenever the
four planes form a transverse quartet. Indeed, all six off-diagonal pairings
are then $\pm1$. Orient the bivectors so that the first row has three
$+1$'s. If the three remaining edge signs agree, $C$ is equivalent up to
sign changes and overall reversal to $J_4-I_4$, where $J_4$ is all ones.
Its eigenvalues are $3,-1,-1,-1$. For mixed signs, permuting the last
three vertices and reversing the pairing and regauging if necessary leaves
one case, with remaining signs $+,+,-$. Its characteristic polynomial
is $(z^2-1)(z^2-5)$, giving $\tr\abs{C}=2+2\sqrt5>6$, impossible.

\begin{planesalgorithmblock}
\begin{mathematicalalgorithm}[Star-free eights without a transverse quartet]\label{alg:starfree-eight}\leavevmode\par\smallskip
\textbf{Input.} The ordered $35$ points of $\mathcal Q$. Join two points
when their reduced wedge pairing is nonzero.
\begin{algorithmsteps}
\item Starting at the empty set, adjoin points in increasing order.
Skip an extension if it contains an independent pairwise orthogonal
triple or four mutually adjacent points (a transverse quartet for realized planes). Record each retained eight-set
and end that branch.
\item In each $\operatorname{GL}(4,2)$-orbit take the lexicographically
least sorted eight-set. Set the diagonal and nonedges of its pairing
matrix to zero. Give every edge of a spanning forest sign $+1$, and
try both signs on every remaining edge.
\item For every three- or four-element subset whose binary points sum to
zero, test whether its columns have a signed sum equal to zero in all eight
rows. Retain the signing if every such subset passes.
\end{algorithmsteps}
\textbf{Return.} The retained signed pairing matrices for each representative.
\end{mathematicalalgorithm}
\end{planesalgorithmblock}

\begin{lemma}[Exhaustiveness of Algorithm~\ref{alg:starfree-eight}]\label{lem:eight-finite}
Every star-free eight of minimum planes without a transverse quartet has,
up to signs and reordering, a returned pairing. Each returned matrix
either has a four-row principal submatrix of trace norm six or agrees,
up to signs and volume reversal, with \eqref{eq:eight-type7-data} or
\eqref{eq:eight-grid-data}.
\end{lemma}
\begin{proof}
Lemma~\ref{lem:binary-lifting} and the binary star criterion justify the
triple test. An integral pairing is $0$ or $\pm1$, so adjacency detects
transversality exactly. Both forbidden patterns persist under extension;
increasing order visits every eligible set once. Elementary basis additions
generate $\operatorname{GL}(4,2)$ and lift integrally. Orienting each tree
outward from a root makes every forest edge positive. Each binary relation
tested in step~3 lifts by Lemma~\ref{lem:binary-lifting} to a signed
vector relation, hence to the required signed sum of columns.

The enumeration gives $27{,}720$ eight-sets in five orbits. Encode the
six exterior coordinates in their stated order, least significant bit
first, and order representatives lexicographically. Scanning edges
lexicographically to build the forest gives the following counts:
\[
\begin{array}{l|rrrrr}
\text{Representative}&1&2&3&4&5\\ \hline
\text{Orbit size}&5040&5040&2520&10080&5040\\
\text{Sign assignments}&256&128&128&128&128\\
\text{Retained}&2&2&2&2&0
\end{array}
\]
For both signings of representative $1$, rows $(1,3,5,7)$, numbered
from zero, have characteristic polynomial $(z^2-1)(z^2-4)$ on their
principal submatrix. Rows $(2,4,5,7)$ give the same polynomial for
representative $4$. In each of these quartets one pairing vanishes, but its
trace norm is still six. The two signings each of
representatives $2$ and $3$ give respectively \eqref{eq:eight-type7-data}
and \eqref{eq:eight-grid-data}, up to bivector signs and volume reversal.
\end{proof}

\subsubsection{Two impossible configurations}
\begin{lemma}\label{lem:eight-exclusions}
The wedge pairings in \eqref{eq:eight-type7-data} and \eqref{eq:eight-grid-data} cannot occur
among eight star-free minimum planes.
\end{lemma}
\begin{proof}

\noindent\emph{First pairing.}
Put $a=x_0$ and $b=x_1$. The first matrix has
$Q(a,b)=Q(a,x_2)=Q(b,x_2)=0$. If $a,b,x_2$ were independent, their planes
would form a star by Lemma~\ref{lem:independent-incident-star}.
Lemma~\ref{lem:binary-lifting} therefore gives a signed relation.
Pairing it with all eight bivectors puts its coefficients in the kernel of
the submatrix of $Q$ on columns $0,1,2$. This kernel is spanned by
$(1,-1,1)$, so $a-b+x_2=0$, or $x_2=b-a$. The remaining vector identities
in \eqref{eq:eight-type7-data} and \eqref{eq:eight-grid-data} follow in the
same way. By bilinearity the displayed data recover the full pairing matrices.

Put $c=x_3,d=x_5,e=x_6$. The data for the first pairing are
\begin{equation}\label{eq:eight-type7-data}
\begin{gathered}
 Q(a,e)=Q(b,c)=Q(b,d)=1,\qquad Q(c,d)=-1,\\
 Q(a,b)=Q(a,c)=Q(a,d)=Q(b,e)=Q(c,e)=Q(d,e)=0,\\
 x_2=b-a,\qquad x_4=a+c,\qquad x_7=d+e.
\end{gathered}
\end{equation}
Put $u=a-b+c$. By \eqref{eq:eight-type7-data}, the quadratic values
of $u+d,u-2d,b+c$ are $-3,3,1$, respectively. Expanding the squared
norms using $\|b-a\|^2=\|a+c\|^2=m$ and applying \eqref{eq:eight-hodge} gives
\[
 15m=2\|u+d\|^2+\|u-2d\|^2+3\|b+c\|^2
 \ge 2V(2\cdot3+3+3)=24V,
\]
contradicting $m\le3V/2$.

\par\noindent\emph{Second pairing.}
Put $a=x_0,b=x_1,c=-x_3,d=x_5$. The data are
\begin{equation}\label{eq:eight-grid-data}
\begin{gathered}
 Q(a,d)=1,\quad Q(b,c)=-1,\quad
 Q(a,b)=Q(a,c)=Q(b,d)=Q(c,d)=0,\\
 x_2=b-a,\quad x_4=a-c,\quad x_6=d-b,\quad x_7=d-c.
\end{gathered}
\end{equation}
In the normalized metric $H$ the eight bivectors have squared length
$t=m/V\le3/2$. Put $p=t-\langle a,d\rangle_H$ and
$q=t-\langle b,c\rangle_H$. The four differences in
\eqref{eq:eight-grid-data} have squared length $t$, so the other four
off-diagonal inner products are $t/2$. Testing $H\pm Q$ on $a-d$ and
$b-c$ gives $p,q\ge1$. The restriction of $Q$ to $a,b,c,d$ is
nonsingular, so $w=a+d-b-c$ is a nonzero integral bivector. The
exterior shortest-vector bound gives
\[
 t\le\|w\|_H^2=2(2t-p-q),\qquad p+q\le 3t/2.
\]

The span $M=\langle a,b,c,d\rangle_{\mathbb Z}$ is unimodular for $Q$,
so $K=M\oplus M^\perp$. The even unimodular complement of signature
$(1,1)$ has an integral null basis $u,v$ with $Q(u,v)=1$. Both have
squared length at least $t$. The trace identity in
\eqref{eq:metric-compatibility} gives $\langle u,v\rangle_H=p-q$.
In the basis $(a,d,b,c,u,v)$, $Q$ is block diagonal for $M\oplus M^\perp$
and $Q^{-1}=Q$. Compatibility gives $H^{-1}=QHQ$ and $\det H=1$.
Jacobi's identity now gives $\det H_M=\det(H^{-1})_{M^\perp}
=\det H_{M^\perp}$, where subscripts denote principal blocks.
On $M$, taking the normalized differences and sums of the pairs $(a,d)$
and $(b,c)$ splits $H_M$ into the scalars $p,q$ and the block
$\left(\begin{smallmatrix}2t-p&t\\t&2t-q\end{smallmatrix}\right)$.
Equating determinants gives
\begin{equation}\label{eq:eight-type11-lower}
 F(p,q):=pq\bigl((2t-p)(2t-q)-t^2\bigr)+(p-q)^2
 =\|u\|_H^2\|v\|_H^2\ge t^2.
\end{equation}

Put $s=(p+q)/2$. Then $2t/3\le1\le s\le3t/4$ and
$\abs{p-q}\le3t/2-2\le t/6$, since $t\le3/2$.
Both diagonal entries and the determinant of the displayed block are
positive. The arithmetic--geometric mean inequality gives $pq\le s^2$
and $(2t-p)(2t-q)\le(2t-s)^2$, hence
\[
 \det H_M\le s^2(t-s)(3t-s).
\]
The three numbers $s/2,s/2,t-s$ sum to $t$, so
$s^2(t-s)\le4t^3/27$. Using $3t-s\le7t/3$ and $t^2\le9/4$
gives $\det H_M\le7t^2/9$. Equation~\eqref{eq:eight-type11-lower}
would now imply
\[
 t^2\le F(p,q)\le(7/9+1/36)t^2=29t^2/36<t^2.\qedhere
\]
\end{proof}

\subsubsection{Star-free subsets of the $D_4$ grid}
\begin{lemma}\label{lem:d4-selected-families}
The lattice $D_4$ has sixteen minimum planes. Its twelve star-free
eight-plane families form one isometry orbit and have exterior rank five; no selected
nine is star-free. It has a star-free seven-plane family of exterior rank six.
\end{lemma}
\begin{proof}
By Lemma~\ref{lem:even-determinant-three}, the minimum bivectors of $D_4$
are wedges of roots $\pm e_i\pm e_j$ with inner product $\pm1$. Such roots
share exactly one coordinate, so their wedge has three entries $\pm1$ on
the pairs of a coordinate triple. On each triple $i<j<k$, start with
$(e_i+e_j)\wedge(e_i+e_k)$. The original wedge and its three single-coordinate sign changes give
four distinct directions on this triple. In coordinates $(12,13,14,23,24,34)$, take the representative
of each pair with first nonzero entry positive and label them
$p_0,\ldots,p_{15}$ in lexicographic order. Three planes
form a star exactly when they occupy distinct rows and distinct columns:
\[
 T=\begin{pmatrix}
 0&6&11&8\\7&1&15&12\\10&14&2&5\\9&13&4&3
 \end{pmatrix},\qquad
 \def\selectedcell{\textcolor[rgb]{.24,.38,.48}{\bullet}}
 \underbrace{\begin{smallmatrix}
 \selectedcell&\selectedcell&\selectedcell&\selectedcell\\\selectedcell&\selectedcell&\selectedcell&\selectedcell\\
 \cdot&\cdot&\cdot&\cdot\\\cdot&\cdot&\cdot&\cdot
 \end{smallmatrix}}_{\text{eight of rank five}},\qquad
 \underbrace{\begin{smallmatrix}
 \selectedcell&\selectedcell&\selectedcell&\selectedcell\\\cdot&\selectedcell&\cdot&\cdot\\
 \cdot&\selectedcell&\cdot&\cdot\\\cdot&\selectedcell&\cdot&\cdot
 \end{smallmatrix}}_{\text{seven of rank six}}.
\]
Here $p_i\wedge p_j=2Q(p_i,p_j)\,e_1\wedge\cdots\wedge e_4$ because
$D_4$ has index two in $\Z^4$. Substitution gives $Q(p_i,p_j)\ne0$
for distinct squares exactly when they share a row or column.
Each bivector has odd coordinate sum, so no signed sum of three can vanish.
Lemma~\ref{lem:binary-lifting}
therefore gives independence for every triple, and
Lemma~\ref{lem:independent-incident-star} proves the grid rule.

Regard each selected square as an edge joining its row and column.
A star is a matching of size three, so by K\H{o}nig's matching
theorem~\cite{Konig1931} a star-free selection is covered by at most two
rows or columns in total. Two rows or two columns cover eight squares,
whereas one of each covers seven. Thus the twelve pairs of rows or of
columns give all star-free eights, and every nine contain a star.

Number the rows and columns from $1$ to $4$. Negating the first
coordinate interchanges row $i$ with column $i$ for every $i$.
Swapping coordinates $1,2$ or $2,3$ and then negating the first
coordinate induces the row permutation $(1\,3\,2\,4)$ or $(1\,4)$,
respectively. These generate $S_4$, so the twelve eights form one
isometry orbit. Each vector in the first two rows annihilates
$(0,-1,1,1,1,0)$ by the coordinate dot product. The first five of
these vectors in row order, restricted to columns $12,13,14,23,34$,
have determinant $4$. Thus the eight have rank five.

The first row and second column form a star-free seven. Six of its
bivectors, in the stated coordinates, satisfy
$\det(p_0,p_1,p_6,p_8,p_{11},p_{13})=-8$, so the seven span the exterior square.
\end{proof}

\subsection{Complete minimum sets when every eight planes contain a star}
\label{app:rankfour-finite-details}\label{sec:proof-algorithms}\label{sec:algorithms}
It remains to consider a rank-four local maximum at which every eight
minimum planes contain a star.
Algorithm~\ref{alg:complete-sets} enumerates candidates for the complete
minimum set. Lemma~\ref{prop:rankfour-configuration-census} proves
exhaustiveness, and the \hyperref[para:metric-outcome]{metric tests}
leave one reference.

\begin{proposition}[Semidefinite lower bounds for selected minimum directions]\label{prop:two-sided-certificate}
Let $S$ be any finite selected family of minimum directions, possibly
partial. Suppose we have
$P,N\succeq0$, nonnegative weights $\beta_c$ on finitely many nonzero
integral null vectors $c\in K$, arbitrary real coefficients $\alpha_x$
for $x\in S$, and a scalar $\gamma$, satisfying
\[
 P+N+\sum_c\beta_ccc^T+\sum_x\alpha_xxx^T+\gamma Q^{-1}=0,
 \qquad \sum_c\beta_c+\sum_x\alpha_x=-1.
\]
Then $t\ge\ell:=\tr(Q(P-N))$. If $S$ spans and a feasible upper matrix $U$ gives
$u=\max_{x\in S}x^TUx<\ell$, then $S$ cannot be the complete minimum
set of a local maximum.
\end{proposition}
\begin{proof}
For $M\pm Q\succeq0$ and $P,N\succeq0$,
\begin{equation}\label{eq:psd-weak-duality}
 \tr(M(P+N))-\tr(Q(P-N))=\tr((M-Q)P)+\tr((M+Q)N)\ge0.
\end{equation}
Pair the first identity of
the proposition with $H$. By \eqref{eq:psd-weak-duality}, $c^THc\ge t$,
$x^THx=t$ and $\tr(HQ^{-1})=0$, its left side is at least
$\ell+t(\sum_c\beta_c+\sum_x\alpha_x)=\ell-t$. This gives the lower
bound. For a complete spanning family Lemma~\ref{lem:unique-config-metric}
gives the upper bound under eutaxy.
\end{proof}

\begin{planesalgorithmblock}
\begin{mathematicalalgorithm}[Candidates for complete minimum sets]\label{alg:complete-sets}\leavevmode\par\smallskip
\textbf{Input.} Basis indices $d=1,2$, the pairing restrictions of
Section~\ref{subsec:finite-basis}, and the index-two glue conditions
\eqref{eq:glue}.
\begin{algorithmsteps}
\item Enumerate the possible wedge-pairing matrices of the six chosen
bivectors, up to reordering, sign changes, and reversal of the pairing.
Apply the trace bound. When $d=2$, enumerate the possible coordinate
lattices for each remaining matrix. Apply the remaining strict and
equality bounds described below.
\item For each surviving pair $(Q,K)$, form the candidate set
\eqref{eq:candidate-cube}. Test each candidate's residue and short relations with the marked six.
Start with the marked six. While the
current layer is nonempty, adjoin each possible new direction. Retain
extensions with pairings in $\{0,\pm1\}$ and a star in every eight.
For each independent six, compute its index in the actual $K$;
require index at most $d$ and its signed pairing class on the surviving
list for that index.
\item Record every retained set of size at least ten. Transport each
whole set to a reference configuration by integral basis exchanges and
symmetries preserving the coordinate lattice and the pairing up to sign.
\end{algorithmsteps}
\textbf{Return.} Reference configurations and a checked transport for
each recorded set.
\end{mathematicalalgorithm}
\end{planesalgorithmblock}

\begin{lemma}[Exhaustiveness of Algorithm~\ref{alg:complete-sets}]\label{prop:rankfour-configuration-census}
Suppose $L$ is a rank-four local maximum and every eight minimum planes
of $L$ contain a star. Algorithm~\ref{alg:complete-sets} records its
complete set of minimum directions in a marked basis.
\end{lemma}
\begin{proof}
Choose six minimum bivectors of maximal index $d$ as in
Section~\ref{subsec:finite-basis}. Their pairing matrix has determinant
$-d^2$. For $d=1$ their coordinate lattice is $K=\Z^6$; for $d=2$
it is $K=\Z^6+\Z(\varepsilon/2)$, where Lemma~\ref{lem:minbasis-index}
and even integrality require $\varepsilon\in\{0,1\}^6$ with
\begin{equation}\label{eq:glue}
 \sum_i\varepsilon_i\in\{5,6\},\qquad
 Q\varepsilon\in2\Z^6,\qquad \varepsilon^TQ\varepsilon\in8\Z.
\end{equation}
In the marked basis the normalized metric $H$ has all six diagonal entries
equal to $t$. Lemma~\ref{lem:trace-reference} applied to $H\pm Q\succeq0$ gives
$\tr\abs{Q}\le\tr H=6t\le9$. The unweighted trace test certifies a violation
by exhibiting a rational symmetric $R$ with operator norm at most one
and $\tr(RQ)>9$.
The strict lower bounds of Proposition~\ref{prop:two-sided-certificate}
exclude pairs $(Q,K)$ for which they force $t>3/2$. A pair forcing
equality in Rankin's bound~\cite[Theorem~5]{Rankin} can occur only for
$D_4$, which contains a star-free eight.

By Lemma~\ref{lem:binary-lifting}, each further minimum direction has
nonzero residue in $K/2K$ distinct from the residues of the marked six, and its three- and
four-term binary relations with them lift to signed vector relations.
For every independent six in an extension, the absolute determinant
in marked coordinates times $d$ gives its index in $K$. Its inherited
glue satisfies \eqref{eq:glue} at index two. Thus every subset of a genuine minimum
family containing the marked six passes all extension tests. A local
maximum has at least ten minimum planes by Proposition~\ref{prop:rankin-direct}.
Recording every retained set from size ten onward therefore includes
the complete family even if it admits a further combinatorial extension.
\end{proof}

\noindent\textbf{Output of Algorithm~\ref{alg:complete-sets}.}
After the unweighted trace test, each index-two pairing class admits at
most one glue up to signed permutations and pairing reversal. The algorithm
discards two pairs that force equality in Rankin's bound.
Of the $25$ surviving pairs $(Q,K)$ in the marked coordinates, $24$ have
$d=1$ and one has $d=2$.
The enumeration records $241{,}807$ marked sets. For each set the exact
check verifies a map $A$ to one of $505$ references $(Q',K',S')$ with
\[
 AK=K',\qquad A^tQ'A=\pm Q,\qquad AS=S'
\]
as unoriented directions in the last equality. Pairing reversal
interchanges $H+Q$ and $H-Q$; the exterior square of a four-dimensional
reflection restores the sign in standard coordinates.

\begin{proposition}[Checking the forced metric]\label{prop:reference-tests}
Let $S$ be a spanning candidate for the complete minimum set and $U$ a
feasible matrix. Suppose a lower bound for $t$ equals
$r=\max_{x\in S}x^TUx$. Any eutactic realization of $S$ has metric $U$.
In particular, $U$ must be compatible, every $x\in S$ must have squared
norm $r$, and every nonzero integral null vector outside $\pm S$ must
have squared norm greater than $r$. At a local maximum the evaluations
$D\mapsto x^TDx$ on the tangent space must span its dual and admit a
strictly positive balance.
\end{proposition}
\begin{proof}
For a eutactic realization, Lemma~\ref{lem:unique-config-metric} gives its
metric as the unique minimizer. Since $U$ attains the lower bound, it is
that metric.

One way to obtain this lower bound is to take
$E_0=\sum_{x\in S}w_xxx^T$ with $w_x\ge0$, $\sum_xw_x=1$ and
$x^TUx=r$ whenever $w_x>0$. If $UE_0U=QE_0Q$, put $R=E_0^{1/2}$.
Then $RUR=\abs{RQR}$, since $RUR\succeq0$ and $(RUR)^2=(RQR)^2$.
For every feasible $M$, Lemma~\ref{lem:trace-reference} applied to
$RMR\pm RQR\succeq0$ gives
\[
 \max_{x\in S}x^TMx\ge\tr(ME_0)\ge\tr(UE_0)=r,
\]
even when $E_0$ is singular.

A single nonzero integral null vector $c\notin\pm S$ with $c^TUc\le r$
therefore excludes $S$. For compatible $U$, the tangent space consists
of the symmetric $D$ satisfying $DQ^{-1}U+UQ^{-1}D=0$.
Corollary~\ref{prop:local-criterion} gives the two tangent conditions.
When the evaluations span, a nonzero tangent with all evaluations
nonnegative rules out a strictly positive balance.
\end{proof}

\phantomsection\label{para:metric-outcome}
\noindent\textbf{Outcome of the metric tests.}
Of the $505$ reference triples $(Q,K,S)$, the strict bounds of
Proposition~\ref{prop:two-sided-certificate} exclude $430$. In twenty of the remaining references, four
selected vectors (the columns of $X$) have wedge-pairing matrix
$C=X^TQX$ with $\tr\abs{C}=6$. Since $X^THX\pm C\succeq0$,
Lemma~\ref{lem:trace-reference} and Rankin's bound give
\[
 6=\tr\abs{C}\le\tr(X^THX)=4t\le6.
\]
Equality forces $D_4$ up to similarity, whose minimum planes contain a
star-free eight by Lemma~\ref{lem:d4-selected-families}.
Proposition~\ref{prop:reference-tests} excludes $54$ more references.
The sole survivor is the index-one reference with ten minimum directions
realized by $A_4$.

For each of the $485$ references outside the twenty quartet cases,
the verified upper matrix satisfies $U\pm Q\succeq0$ and
$x^TUx\le u<3/2$ on every selected direction. With all directions as
columns of $X$, congruence and Lemma~\ref{lem:trace-reference} give,
for every nonzero real diagonal matrix $D$,
\[
 \tr\abs{DX^TQXD}\le\tr(DX^TUXD)\le u\tr(D^2)<\tfrac32\tr(D^2).
\]
Zero diagonal entries allow arbitrary subsets, so even weighted subsets
cannot reach the Rankin threshold on these references.

\subsection{Recovering the lattice up to similarity and duality}\label{sec:narain-transport}
Put $E=\mathbb Z^4$ with its standard orientation and $K=\Lambda^2E$.
In the basis $(e_{12},e_{13},e_{14},e_{34},-e_{24},e_{23})$, the wedge
pairing $Q$ has matrix $J=\left(\begin{smallmatrix}0&I_3\\I_3&0\end{smallmatrix}\right)$.
Thus $K$ is even unimodular of signature $(3,3)$, with quadratic value $s(x)=x^TJx/2$.

\begin{proposition}[Integral isometries of the wedge pairing]\label{prop:integral-exterior-equivalence}
The full integral isometry group is
\begin{equation}\label{eq:integral-exterior-group}
 O(K,Q)=\{\,\pm\Lambda^2g,\ \pm J\Lambda^2g:
 g\in\operatorname{SL}_4(\mathbb Z)\,\}.
\end{equation}
\end{proposition}
\begin{proof}
Shioda's integral exterior-square result, as stated
in~\cite[\S4.1, before Theorem~4.1.1]{LiZou2026}, gives
$\Lambda^2\operatorname{SL}_4(\mathbb Z)$ index two in $\operatorname{SO}(K,Q)$,
with $-I_6$ outside its image. Thus the determinant-one isometries are
$\pm\Lambda^2 g$. Since $J$ exchanges three pairs of basis vectors,
it is an integral isometry of determinant $-1$, giving the other coset.
\end{proof}

For a positive Gram matrix $G$ normalized by $\det G=1$, let $H=\Lambda^2G$.
The maps $\pm\Lambda^2g$ act on $H$ by the basis change $G\mapsto g^TGg$.
Since $G$ preserves volume, $H^TJH=J$, whence
\begin{equation}\label{eq:narain-duality}
 J^T(\Lambda^2G)J=H^{-1}=\Lambda^2(G^{-1}).
\end{equation}
Thus $J$ acts by duality. To see that $H$ determines $G$, suppose
$\Lambda^2G_1=\Lambda^2G_2$. Then $P=G_1^{-1/2}G_2G_1^{-1/2}>0$
satisfies $\Lambda^2P=I$. All pairwise products of its four positive
eigenvalues are one, so $P=I$ and $G_1=G_2$. The exterior metric and
integral pairing therefore recover the oriented covolume-one lattice up
to duality. For an integral map reversing the pairing, first reverse
the orientation of one rank-four lattice.

\paragraph{Acknowledgments.} I thank Henry Cohn for posing the question about
spectral gaps of Narain lattices that led me to this problem, and for his
mentorship when I began that work in 2020.

\section*{Computational resources and use of AI tools}
\addcontentsline{toc}{section}{Computational resources and use of AI tools}
The code, exact witnesses and reproduction instructions accompany the paper as
arXiv ancillary files. Version~1.0.0 is archived on Zenodo~\cite{MehraMinimumPlanesComputations}.
The source is also available on
\href{https://github.com/pihart/minimum-planes-computations/tree/v1.0.0}{GitHub}.
The script \texttt{check\_eight.py} runs the enumeration in
Algorithm~\ref{alg:starfree-eight} and verifies the bivector ranks in the $D_4$ grid examples.
The script \texttt{check\_quartets.py} checks the final classification of
$505$ reference cases. Strict gaps exclude $430$ cases, quartets forcing
Rankin equality exclude $20$, and further metric tests exclude $54$.
The remaining case gives the ten minimum planes of $A_4$.

The author developed the mathematical arguments largely independently of
AI assistance. AI tools helped implement and verify the finite enumerations
and metric bounds, review and simplify proofs, locate and check references,
and draft and revise the exposition.

\bibliographystyle{plainurl}\bibliography{references}
\end{document}